\documentclass[12pt,twoside,reqno]{article}

\usepackage{amsmath}
\usepackage{latexsym}
\usepackage{amsmath,amsthm}
\usepackage{amsfonts}
\usepackage{amssymb}
\usepackage{amsfonts,euscript}
\usepackage{graphicx}
\newtheorem{theorem}{\bf Theorem}[section]

\newtheorem{lemma}{\bf Lemma}[section]

\newcommand{\beq}{\begin{equation}}
\newcommand{\eeq}{\end{equation}}
\newcommand{\beqn}{\begin{eqnarray}}
\newcommand{\eeqn}{\end{eqnarray}}
\newcommand{\bear}{\begin{array}}
\newcommand{\eear}{\end{array}}
\newcommand{\beit}{\begin{itemize}}
\newcommand{\eeit}{\end{itemize}}
\newcommand{\beqno}{\begin{eqnarray*}}
\newcommand{\eeqno}{\end{eqnarray*}}

\let\theta\vartheta
\def\eqref #1{(\ref{#1})}

\allowdisplaybreaks

\numberwithin{equation}{section}
\begin{document}
\thispagestyle{plain}
\title{Design of boundary stabilizers for the non-autonomous cubic semilinear heat  equation, driven by a multiplicative noise}
\maketitle

\begin{center}
\vskip-1cm

{\large\sc Ionut Munteanu}\\
{\normalsize e-mail: {\tt ionut.munteanu@uaic.ro}}\\

{\small Alexandru Ioan Cuza University, Department of Mathematics\\
 and Octav Mayer Institute of Mathematics (Romanian Academy)\\ 
  700506 Ia\c{s}i, Romania}
\end{center}

\begin{abstract}Here we design  boundary feedback stabilizers to unbounded trajectories, for semi-linear stochastic heat equation with  cubic non-linearity.   The feedback controller is linear, given in a simple explicit form and involves only the eigenfunctions of the Laplace operator. It is supported in a given open subset of the boundary of the domain. Via a rescaling argument, we transform the stochastic equation into a random deterministic one. Then, the simple form of the feedback, we propose here, allows to write the solution, of the random equation, in a mild formulation via a kernel. Appealing to a fixed point argument the existence \& stabilization result is proved.  

\noindent \textbf{\textit{MSC 2010 Classification:}}

\noindent \textbf{\textit{Keywords:}}
\end{abstract}
\section{Presentation of the model}
Let $\mathcal{O}\subset\mathbb{R}^2$ be a bounded domain, with its smooth boundary $\partial\mathcal{O}$ split in two parts as $\partial\mathcal{O}=\Gamma_1\cup\Gamma_2$, such that  $\Gamma_1$ has non-zero surface measure. We consider the following boundary controlled semi-linear heat equation, with cubic non-linearity, driven by a multiplicative noise
\begin{equation}\label{2m2}\left\{\begin{array}{l}\begin{aligned} dY(t,x)=(\Delta Y(t,x)+cY(t,x)&+f(t,x,Y(t,x)))dt
+\theta Y(t,x)dW(t),\\&  \text{for } t>0,\ x\in \mathcal{O},\end{aligned}\\
\\
\frac{\partial}{\partial\mathbf{n}}Y(t,x)=u(t,x), \text{ on }t\geq0,\ x\in\Gamma_1,\\
\\
\frac{\partial}{\partial\mathbf{n}}Y(t,x)=0,\ \text{ on } t\geq0,\ x\in \Gamma_2,\\
\\
Y(0,x)=y_o(x),\ x\in \mathcal{O}.\end{array}\right.\ \end{equation}
Here, $dW$ denotes a Gaussian time noise, that is usually understood as the distribution derivative of the Brownian sheet $W(t)$ on a probability space $(\Omega, \mathcal{F},\mathbb{P})$ with normal filtration $(\mathcal{F}_t)_{t\geq0}$. $c$ and $\theta$ are some positive constants. $f$ is a cubic polynomial with time-space coefficients, of the form
$$f(t,x,y)=a_2(t,x)y^2+a_3(t,x)y^3.$$ 
   On the functions $a_i,\ i=2,3,$ we assume that: there exist $C_a>0$ and 
$$ 0\leq m^i_1\leq m^i_2\leq ...\leq m^i_{S_i},$$ for some $S_i\in\mathbb{N}$, $i=2,3,$ for which
\begin{equation}\label{a}(H_0)\ \ \  \ \sup_{x\in\mathcal{O}}|a_i(t,x)|\leq C_a\left(\sum_{k=1}^{S_i}t^{m^i_k}+1\right),\ \forall t\geq 0,\ i=2,3.\end{equation} Moreover, we assume that $m^2_{S_2}$, $m^3_{S_3}$ and $\theta$ are such that:  
\begin{equation}\label{gogol32}(H_1)\ \ \ \ \ \  \ \ \  \ \ \ \ \ \ \ \ \, \frac{1}{2}\theta^2-m_S-\frac{1}{100}=\theta_1>0,  \ \ \ \ \ \ \ \ \ \ \ \ \  \ \ \ \ \ \ \ \ \end{equation}where  $m_S:=\max\left\{m^2_{S_2},m^3_{S_3}\right\}$. 

\noindent (We remark that, when $a_2\equiv 0$, then we stumble exactly on the non-autonomous Chafee-Infante equation.)

 \noindent Finally,  $\mathbf{n}$ stands for the outward unit normal to the boundary $\partial\mathcal{O}$, and $u$ is the control. The initial data $y_o$ is $\mathcal{F}_0$-adapted.

The aim of the present paper is to find  a feedback law $u$ such that, once inserted into the equation (\ref{2m2}), the corresponding solution of the closed-loop equation (\ref{2m2}) satisfies
$$ e^{\alpha t}\int_\mathcal{O}Y^2(t,x)dx<const.,\ \forall t\geq0, \mathbb{P}-\text{a.s.},$$for a prescribed positive constant $\alpha$, provided that the initial data $y_o$ is small enough in the $L^2-$norm (that is the main result stated in Theorem \ref{t1} below). Note that this  is an almost sure path-wise local boundary stabilization type result. Besides this, since the coefficients are time-dependent, our considerations are related, in fact, to the problem of stabilization to trajectories (i.e., non-steady states). In the existing literature there are only few results on this problem. Regarding the internal stabilization we refer to \cite{btraj}, while for the the boundary case we cite \cite{rodrigues1,rodrigues2,ion8}.  In any case,  the time-dependent coefficients are assumed to be bounded, while here, we let them explode when $t$ goes to infinity. This, together with the   noise perturbations, makes our task a lot more difficult. Note that even the well-posedness is not known for our example. Anyway, the simple form of the controller, which we shall introduce below, allows us to write the equations in an integral formulation, via a kernel. Then, via a fixed point argument and a proper choose of the spaces, the three raised problems, i.e., existence, uniqueness and stabilization, will be solved. 

It is worth to mention that the work \cite{caraballo} studies the effect of noise on the Chafee-Infante equation, and the conclusions there state that a single multiplicative Ito noise, of sufficient intensity, will stabilize the origin of the system. However, we remark that the coefficients there are assumed to be bounded, and then, the "sufficient intensity" of the noise is related to their bounds. While here,  due to the unboundedness of the coefficients, those arguments cannot be applied. Anyway, the presence of the noise is mandatory. This can be seen from the imposed hypothesis \eqref{gogol32}. But, even if the level of the noise, $\theta$, is large, it cannot ensure the stability of the system. A boundary stabilizer is needed. In conclusion, the result of this paper is  first in this general framework.

 It is clear that, due to the general form of the nonlinearity $f$, the results presented here can be applied not only to the Chafee-Infante equation, since, many examples of cubic semi-linear equation arise from biology, chemistry or physics, such as the FitzHugh-Nagumo model \cite{fitzhugh}(in neuroscience) or the Fischer-Kolmogorov model \cite{graham} (in evolution of population dynamics).

 The method to design the feedback controller  $u$ relies on the ideas in \cite{ion1}, where a proportional type law was proposed to stabilize, in mean, the stochastic heat equation. We emphasize that, unlike to the  equation in \cite{ion1}, which is linear and evolves in a bounded interval, now we deal with a nonlinear one of order three, evolving in the 2-D domain $\mathcal{O}$.  In order to overcome this  complexity, we further develop the ideas in \cite{ion1}. Roughly saying, we  design a similar feedback as in \cite{ion1}: linear, of finite-dimensional structure, given in a very simple form, being easy to manipulate from the computational point of view, involving only the eigenfunctions of the Neumann-Laplace operator (see relations \eqref{e1}-\eqref{u} below). Then, we plug it into the equations, and show that it achieves the stability by using the estimates on the magnitude of the controller  and a fixed point argument in a properly chosen space (see Theorem \ref{t1} below). The idea to use fixed point arguments in order to show the stability of deterministic or stochastic equations has been previously used in papers like \cite{fix2,fix1}. Proportional type feedback, similar to that one we  design here, has its origins in the works  \cite{barbui,ion2}, while in the papers \cite{ion3,ion4,ion5,ion6,ion7}, it has been used to stabilize other  important parabolic-type equations, such as the Navier-Stokes equations (also with delays), the Magnetohydrodynaimc equations, or the phase field equations.  Besides the method of proportional-type controllers, the backstepping technique has been developed with lots of important results. Even if, at a first glance, the two methods seem to be very similar, conceptually they are totally different. For more details, we refer to the works  \cite{bal,k1}, while in \cite{k2} it is provided also a stabilization result for the stochastic Burgers equation.

\section{The random equation}

There is a well-known trick, by now, on how to avoid to deal with stochastic equations. Namely, to equivalently rewrite  them as random deterministic ones via a rescaling argument. This is explained in full details in the work \cite{barbur}. To this end, in (\ref{2m2}), let us consider the transformation
\begin{equation}Y(t)=\Gamma (t)y(t),\ t\in[0,\infty),\end{equation}where $\Gamma(t):L^2(\mathcal{O})\rightarrow L^2(\mathcal{O})$ is the linear continuous operator defined by the equations
$$d\Gamma(t)=\theta\Gamma(t)dW(t),\ t\geq0,\ \Gamma(0)=1,$$that can be equivalently expressed as
\begin{equation}\Gamma(t)=e^{\theta W(t)-\frac{t}{2}\theta^2},\ t\geq0.\end{equation}
Frequently below we shall use the obvious inequality
$$e^{-at}\leq t^{-a},\ \forall t>0,a\geq0.$$
 By the law of the iterated logarithm and arguing similarly as in Lemma 3.4 in \cite{barbup}, it follows that there exists a constant $C_\Gamma>0$ such that 
\begin{equation}\label{iv10}\Gamma(t)= e^{\theta W(t)-\theta_1t}e^{-(m_S+\frac{1}{100}) t}\leq C_\Gamma e^{-(m_S+\frac{1}{100}) t} ,\  \forall t>0, \mathbb{P}-\text{a.s.},\end{equation}where we have used that $\frac{1}{2}\theta^2=m_S+\frac{1}{100}+\theta_1$ assumed in \eqref{gogol32}.
Then, taking advantage of  \eqref{a}, we deduce that, for $i=2,3$, we have
\begin{equation}\label{gogolv10}\begin{aligned}\Gamma(t)\sup_{x\in\mathcal{O}}|a_i(t,x)|&\leq C_\Gamma C_a\left(\sum_{k=1}^{S_i}t^{m^i_k} e^{-(m_S+\frac{1}{100})t}+e^{-\left(m_S+\frac{1}{100}\right)t}\right)\\&
(\text{since $0\leq m^i_1\leq m^i_2\leq ...\leq m^i_{S_i}\leq m_S$})\\&
\leq C \left(\sum_{k=1}^{S_i}t^{m^i_k}e^{-(m^i_k+\frac{1}{100})t}+e^{-\frac{1}{100}t}\right)\\&
\leq C\left(\sum_{k=1}^{S_i}t^{m^i_k}t^{-(m^i_k+\frac{1}{100})}+t^{-\frac{1}{100}}\right)\\&
\leq C t^{-\frac{1}{100}},\ \forall t>0.\end{aligned}\end{equation}
Next, applying  It\"o's formula in (\ref{2m2}) , we obtain that $y$ satisfies the following  random partial differential equation

\begin{equation}\label{m2}\left\{\begin{array}{l}\begin{aligned} \partial_ty(t,x)=\Delta y(t,x)+cy(t,x)&+\Gamma^{-1}(t)f(t,x,\Gamma(t)y(t,x)),\\&
 \text{for } t>0,\ x\in \mathcal{O},\end{aligned}\\
\\
\frac{\partial}{\partial\mathbf{n}}y(t,x)=u(t,x), \text{ on }t\geq0,\ x\in\Gamma_1,\\
\\
\frac{\partial}{\partial\mathbf{n}}y(t,x)=0,\ \text{ on } t\geq0,\ x\in \Gamma_2,\\
\\
y(0,x)=y_o(x),\ x\in \mathcal{O}.\end{array}\right.\ \end{equation}

\section{The boundary feedback stabilizer and the main result of the work}

Let $(X,\|\cdot\|_X)$ stand for some normed space. We set $C_b([0,\infty),X)$ for the space of all continues $X$-valued functions, that are $\|\cdot\|_X-$bounded on $[0,\infty)$. Next, we denote by $L^p$, $1\leq p\leq \infty$, the  Lebesgue space $L^p(\mathcal{O})$ consisting of all power $p$ integrable functions, endowed with the standard norm $|\cdot|_p$; we denote by $W^{s,p},\ s\in (0,1)$ the corresponding fractional Sobolev space, i.e.
$$W^{s,p}(\mathcal{O}):=\left\{y\in L^p(\mathcal{O}):\ \frac{|y(x)-y(\xi)|}{|x-\xi|^{\frac{2}{p}+s}}\in L^p(\mathcal{O}\times\mathcal{O})\right\},$$ which is an intermediary Banach space between $L^p$ and $W^{1,p}(\mathcal{O})$, endowed with the natural norm
$$\|y\|_{s,p}:=\left(\int_\mathcal{O}|y|^pdx+\int_{\mathcal{O}}\int_\mathcal{O}\frac{|y(x)-y(\xi)|^p}{|x-\xi|^{2+sp}}dxd\xi\right)^\frac{1}{p}.$$ For the particular case $p=2$, we set $H^s:= W^{s,2}$ and
$$\|\cdot\|_s:=\|\cdot\|_{s,2}.$$ By interpolation, one can extend the definition of $H^s$ for each $s>0$. We set $H_0^1(\mathcal{O})$ for the completion of the $C_0^\infty(\mathcal{O})$ ( the set of $C^\infty$-compact supported functions in $\mathcal{O}$) in the $H^1$-norm.

It is well known the following fractional Sobolev embedding (for details see \cite{brezis})
\begin{equation}\label{toto1}|y|_4\leq C\|f\|_\frac{1}{2},\ \forall f\in H^{\frac{1}{2}}(\mathcal{O}),\end{equation}where $C=C(\mathcal{O}) $ is some positive constant, depending on the domain $\mathcal{O}$.

Finally, we set $\left<\cdot,\cdot\right>$ for the natural scalar product in $L^2$; and $\left<\cdot,\cdot\right>_N$, the euclidean scalar product in $\mathbb{R}^N,\ N\in\mathbb{N}$. We shall denote by $C$ different constants that may change from line to line, though we keep denote them by the same letter $C$, for the sake of the simplicity of the writing.

 Let us denote by 
$$\mathcal{A}y=-(\Delta y+cy),\ \forall y\in \mathcal{D}(\mathcal{A}),$$
$$\mathcal{D}(\mathcal{A})=\left\{y\in H^2(\mathcal{O}):\, \frac{\partial}{\partial\mathbf{n}}y=0 \text{ on }\partial\mathcal{O}\right\}.$$ Here, $-\Delta$ is  the Neumann-Laplace operator  on $\mathcal{O}.$ It is well known that it has a discrete spectrum, i.e., it has a countable set of semi-simple non-negative eigenvalues $\left\{\lambda_j\right\}_{j=1}^\infty$ with $\lambda_1=0$. We assume that the eigenvalues set is arranged as an increasing sequence with $\lambda_j\rightarrow\infty$ when $j\rightarrow\infty$. We denote by $\left\{\varphi_j\right\}_{j=1}^\infty$, the corresponding eigenfunctions, which form an orthonormal basis in $L^2$.  More precisely, we have
$$-\Delta\varphi_j=\lambda_j\varphi_j \text{ in }\mathcal{O}, \frac{\partial}{\partial\mathbf{n}}\varphi_j=0 \text{ on }\partial\mathcal{O},\ \forall j=1,2,3,...,$$and
$$\left<\varphi_i,\varphi_j\right>=\delta_{ij},\ \forall i,j=1,2,3,...,$$$\delta_{ij}$ being the Kronecker symbol.  Besides this, by the Parseval's identity, we have for a function $y\in L^2$, the following decomposition
$$y=\sum_{j=1}^\infty \left<y,\varphi_j\right>\varphi_j$$and
$$|y|_2=\left(\sum_{j=1}^\infty |\left<y,\varphi_j\right>|^2\right)^\frac{1}{2},$$where $\left<y,\varphi_j\right>,\ j=1,2,3,...,$ are called the (Fourier) modes of $y$. Moreover, since $\mathcal{O}$ is bounded with smooth boundary, it is also known that the norm $\|\cdot\|_\alpha$ is equivalent with $|(-\Delta)^\frac{\alpha}{2}\cdot|_2,$ for all $\alpha>0.$ Thus, one can find some constants $C_1,C_2>0$ such that
\begin{equation}\label{toto3}C_1\left(\sum_{j=1}^\infty \lambda_j^\frac{1}{2}|\left<y,\varphi_j\right>|^2 \right)^\frac{1}{2}\leq \|y\|_\frac{1}{2}\leq C_2\left(\sum_{j=1}^\infty \lambda_j^\frac{1}{2}|\left<y,\varphi_j\right>|^2 \right)^\frac{1}{2},\ \forall y\in W^{\frac{1}{2},2}.\end{equation}

In this work, we shall assume that the eigenvalues system $\left\{\lambda_j\right\}_j$ obeys
\begin{equation}\label{h1} (H_1)\ \ \ \ \ \ \ \  \ \ \ \ \ \ \ \ \ \ \  \sum_{j=2}^\infty\frac{1}{\lambda^\frac{5}{3}_j}<\infty. \ \ \ \ \ \ \ \ \ \ \ \ \end{equation}
In the Appendix below, we shall verify that, when $\mathcal{O}$ is a square, then $(H_1)$ holds true. But one can easily find many more examples of domains $\mathcal{O}$ for which assumption $(H_1)$ is full-filed. 

We go on and recall  the well-known  $L^\infty$-bounds of the Laplace eigenfunctions
\begin{equation}\label{toto2}|\varphi_j|_\infty\leq C\lambda_j^\frac{1}{4},\ \forall j=2,3,...,\end{equation} that hold true without making any geometric assumption on the domain $\mathcal{O}\subset\mathbb{R}^2$. We are also aware of Tataru's trace estimates 
\begin{equation}\label{tataru}\left\|\varphi_j\right\|_{L^2(\partial\mathcal{O})}\leq C\lambda_j^\frac{1}{6},\ j=2,3....\end{equation} 

For latter purpose, let us show that, for an arbitrary constant $C>0$ and a sufficiently large $M>0$, we have that
\begin{equation}\label{i6}\begin{aligned}\sum_{j=M}^\infty e^{(-2\lambda_j+C) t}\varphi^2_j(\xi)\leq C_0\frac{1}{t},\ \forall t>0,\ \xi\in \mathcal{O}.\end{aligned}\end{equation} Here,  $C_0>0$ is some constant. Indeed, since $M$ is large enough and $\lim_{j\rightarrow \infty} \lambda_j=\infty$, we have that
$$-2\lambda_j+C\leq -\lambda_j,\ \forall j\geq M.$$Thus
$$\sum_{j=M}^\infty e^{(-2\lambda_j+C) t}\varphi^2_j(\xi)\leq\sum_{j=M}^\infty e^{-\lambda_j t}\varphi^2_j(\xi).$$ The latter term is the the rest of order $M$ of the well-known Neumann heat kernel, which is known to be less or equal of some constant times $\frac{1}{t}$. From this, our claim \eqref{i6} follows immediately.

Now, let us come back to  the above defined operator $\mathcal{A}$. It is clear that it has as-well discrete semi-simple spectrum, namely
$$\mu_j:=\lambda_j-c,\ j=1,2,3... $$with the corresponding eigenfunctions $\left\{\varphi_j\right\}_{j=1}^\infty.$ 

Let some $\rho>1$. Then, there exists $N\in\mathbb{N}$ such that 
$$\mu_j\leq \rho,\ j=1,2,...,N \text{ and }\mu_j> \rho,\ \forall j\geq N+1.$$ The first $N$ eigenvalues are usually called \textit{the unstable eigenvalues.}

It is obvious that, given any prescribed $\alpha>0$, if we take $\rho$ and $N$ large enough, we may suppose that  the following relations hold true:
\begin{equation}\label{oo}(H_3) \ \ \ \ \ \begin{array}{l}
1)\ -\rho+2\alpha+\frac{3}{4}+\frac{1}{4}\lambda_i-\frac{1}{100}\leq 0;\\
\\
2)\ -\rho+3\alpha+\frac{7}{12}+\frac{1}{4}\lambda_i-\frac{1}{100}\leq 0;\\
\\
3)\ -\rho+2\alpha+\frac{1}{4}(\lambda_i+\lambda_j)+\frac{1}{2}-\frac{1}{100}\leq 0;\\
\\
4)\ -\rho+3\alpha+\frac{5}{12}+\frac{1}{4}(\lambda_i+\lambda_j)-\frac{1}{100}\leq 0;\end{array}\end{equation}for all $i,j=1,2,...,N$.

Although it would be possible to treat the case of semi-simple unstable eigenvalues following
\cite{ion2}, for the sake of simplicity, we assume that
\begin{equation}\label{h2}(H_4)\ \ \ \ \text{ The first $N$ eigenvalues $\mu_j,\ j=1,2,...,N,$ are simple},\end{equation}i.e., we have
$$\mu_1<\mu_2<...<\mu_N.$$

Now, since we are set with the theoretical results and the hypotheses of the paper, we may proceed to apply the approach from \cite{ion1},\cite{ion2}. Firstly, in order to lift the boundary control into the equations (to obtain an internal control-type problem), we introduce  the so-called  Neumann operator as: given $g\in L^2(\Gamma_1)$ and $\gamma>0$, we denote by $D_{\gamma}g:=y$, the solution to  the equation
\begin{equation}\label{e41}\begin{aligned}&- \Delta y(x)-cy(x)-2\sum_{i=1}^N\mu_i\left<y,\varphi_i\right>\varphi_i(x)+\gamma y(x)=0 \\&
\text{ for }x\in\mathcal{O};\  \frac{\partial}{\partial\mathbf{n}}y(x)=g \text{ on }\Gamma_1 \text{ and }\frac{\partial}{\partial\mathbf{n}}y(x)=0 \text{ on }\Gamma_2.\end{aligned}\end{equation} For $\gamma$ large enough, equation \eqref{e41} has a unique solution, defining so the map $D_\gamma\in L(L^2(\Gamma_1),H^\frac{1}{2})$(for further details check e.g. \cite[p. 6]{10}). Also, appealing to  Green formula (see the computations  in \cite[Eqs. (4.1)-(4.2)]{ion2}), we deduce that
\begin{equation}\label{e7}\left<D_\gamma g,\varphi_j\right>=\left\{\begin{array}{l}-\frac{1}{\gamma-\mu_j}\int_{\Gamma_1} g \varphi_jd\sigma,\ j=1,2,...,N,\\
\\-\frac{1}{\gamma+\mu_j}\int_{\Gamma_1} g\varphi_jd\sigma,\ j>N.
\end{array}\right.\ \end{equation}Here, $d\sigma$ is the surface measure on $\Gamma_1$.

 Next, we choose 
$$\gamma_N>\gamma_{N-1}> \dots >\gamma_1>\rho,$$ $N$ constants, large enough, such that equation \eqref{e41} is well-posed for each of them, and denote by $D_{\gamma_i},\ i=1,2,...,N$, the corresponding Neumann maps.

Following the ideas in \cite{ion2}, we denote by $\textbf{B}$ the Gram matrix of the system $\left\{\varphi_i|_{\Gamma_1}\right\}_{i=1}^N$ in the Hilbert space $L^2(\Gamma_1)$, with the standard scalar product  
$$\left<g,h\right>_0:=\int_{\Gamma_1}g(x)h(x)d\sigma.$$More precisely,
  \begin{equation}\label{e9}\textbf{B}:=\left(\begin{array}{cccc}\left<\varphi_1,\varphi_1\right>_0& \left<\varphi_1,\varphi_2\right>_0& \dots &\left<\varphi_1,\varphi_N\right>_0\\
\left<\varphi_2,\varphi_1\right>_0& \left<\varphi_2,\varphi_2\right>_0& \dots &\left<\varphi_2,\varphi_N\right>_0\\
\dots &\dots &\dots & \dots \\
\left<\varphi_N,\varphi_1\right>_0&\left<\varphi_N,\varphi_2\right>_0& \dots &\left<\varphi_N,\varphi_N\right>_0\end{array}\right).\end{equation}Further, we introduce the matrices
\begin{equation}\Lambda_{\gamma_k}:=\left(\begin{array}{cccc}\frac{1}{\gamma_k-\mu_1}&0& \dots &0\\
0&\frac{1}{\gamma_k-\mu_2}& \dots &0
\\\dots &\dots & \dots & \dots \\
0&0& \dots &\frac{1}{\gamma_k-\mu_N} \end{array}\right),\ k=1,...,N,
\end{equation} 
 \begin{equation}\label{T}T:=\left(\begin{array}{cccc}\frac{1}{\gamma_1-\mu_1}\varphi_1|_{\Gamma_1} & \frac{1}{\gamma_1-\mu_2}\varphi_2|_{\Gamma_1} & \dots & \frac{1}{\gamma_1-\mu_N}\varphi_N|_{\Gamma_1}\\
\frac{1}{\gamma_2-\mu_1}\varphi_1|_{\Gamma_1} & \frac{1}{\gamma_2-\mu_2}\varphi_2|_{\Gamma_1} &\dots& \frac{1}{\gamma_2-\mu_N}\varphi_N|_{\Gamma_1}\\
\dots&\dots&\dots&\dots\\
\frac{1}{\gamma_N-\mu_1}\varphi_1|_{\Gamma_1} & \frac{1}{\gamma_N-\mu_2}\varphi_2|_{\Gamma_1} &\dots& \frac{1}{\gamma_N-\mu_N}\varphi_N|_{\Gamma_1}\end{array}\right),\end{equation}and \begin{equation}\label{AA}A=(B_1+B_2+\dots+B_N)^{-1},\end{equation}where  
\begin{equation}\label{bo}B_k:=\Lambda_{\gamma_k}\textbf{B}\Lambda_{\gamma_k},\ k=1,...,N.\end{equation}
We recall the Appendix in \cite{ion2} where it is shown that  the sum $B_1+B_2+\dots +B_N$ is an invertible matrix, and consequently, the matrix $A$ is well-defined.

Now, let us introduce the feedback laws:

\begin{equation}\label{e1}u_k(y)(t,x)=\left<A\left(\begin{array}{c}\left<y(t),\varphi_1\right>\\
\left<y(t),\varphi_2\right>\\\dots \\\left<y(t),\varphi_N\right>\end{array}\right),\left(\begin{array}{c}\frac{1}{\gamma_k-\mu_1}\varphi_1(x)\\ \frac{1}{\gamma_k-\mu_2}\varphi_2(x)\\ \dots \\\frac{1}{\gamma_k-\mu_N}\varphi_N(x)\end{array}\right)\right>_N,\end{equation} for $t\geq0,\ x\in \Gamma_1$, and $k=1,2,...,N$. Then, define $u=u(y)$ as
\begin{equation}\label{toto10}u=u_1+u_2+\dots +u_N,\end{equation}which, in a condensed form, can be written as 
 \begin{equation}\label{u}u= \left<T\ A\left(\begin{array}{c}\left<y(t),\varphi_1\right>\\
\left<y(t),\varphi_2\right>\\\dots \\\left<y(t),\varphi_N\right>\end{array}\right),\left(\begin{array}{c}1\\ 1\\ \dots \\1\end{array}\right)\right>_N.\end{equation} We claim that, once inserted this feedback form $u$ into the equation \eqref{m2} it yields the local exponential asymptotic stability of the corresponding closed-loop system \eqref{m2}. More exactly, we will show that
\begin{theorem}\label{t1} Let $\eta>0$ be sufficiently small. Under $(H_0)$-$(H_4)$, for each $y_o\in L^2$ with $|y_0|_2<\eta$, once plugged the feedback law $u$, given by \eqref{u}, into the equation \eqref{m2}, there exists a unique solution $y$ to the  closed-loop equation \eqref{m2}, which  belongs to the space 
$$\mathcal{Y}:=\left\{y\in C_b([0,\infty),H^\frac{1}{2}(\mathcal{O})):\ \sup_{t>0}\left[e^{\alpha t}(|y(t)|_2+t^{\frac{1}{12}}\|y(t)\|_\frac{1}{2})\right]<\infty\right\}.$$

Consequently, $Y(t)=\Gamma(t)y(t)$ is the unique solution of the stochastic cubic equation \eqref{2m2}, which satisfies
$$e^{\alpha t}\int_\mathcal{O}Y^2(t,x)dx <const.,\ \forall t\geq0,\ \mathbb{P}-\text{a.s.}.$$

\end{theorem}

\section{Proof of the main result}

In order to ease our problem, we shall equivalently rewrite  equation \eqref{m2} as an internal control-type problem, by using similar arguments as in \cite[Eqs. (17)-(19)]{ion1}. We arrive to:
\begin{equation}\label{v1}\begin{aligned}\partial_t y(t)=&-\mathcal{A}y(t)+\sum_{i=1}^N(\mathcal{A}+\gamma_i)D_{\gamma_i}u_i(y(t))\\&-2\sum_{i,j=1}^N\mu_j\left<D_{\gamma_i}u_i(y(t)),\varphi_j\right>\varphi_j \\&
+\Gamma(t)a_2(t)y^2(t)+\left(\Gamma(t)\right)^2a_3(t)y^3(t);\ y(0)=y_o.\end{aligned}\end{equation}

The following result is related to the linear operator which governs equation (\ref{v1}), i.e.
$$\mathbb{A}y:=-\mathcal{A}y(t)+\sum_{i=1}^N(\mathcal{A}+\gamma_i)D_{\gamma_i}u_i(y(t))-2\sum_{i,j=1}^N\mu_j\left<D_{\gamma_i}u_i(y(t)),\varphi_j\right>\varphi_j,$$ $\forall y\in\mathcal{D}(\mathbb{A})=\mathcal{D}(\mathcal{A}).$ It says that the semigroup generated by it can be written in a mild formulation via a kernel $p$, as
$$(e^{t\mathbb{A}}y_o)(x)=\int_\mathcal{O}p(t,x,\xi)y_o(\xi)d\xi,\ t\geq0,\ x\in\mathcal{O}.$$ Its proof is given in the Appendix.
\begin{lemma}\label{l1}The solution $z$ of
\begin{equation}\label{v2}\begin{aligned}&\partial_t z(t)=-\mathcal{A}z(t)+\sum_{i=1}^N(\mathcal{A}+\gamma_i)D_{\gamma_i}u_i(z(t))-2\sum_{i,j=1}^N\mu_j\left<D_{\gamma_i}u_i(z(t)),\varphi_j\right>\varphi_j;\\&
 z(0)=z_o, \end{aligned}\end{equation} can be written in a mild formulation as
$$z(t,x)=\int_\mathcal{O}p(t,x,\xi)z_o(\xi)d\xi,$$where
\begin{equation}\label{v3}p(t,x,\xi):=p_1(t,x,\xi)+p_2(t,x,\xi)+p_3(t,x,\xi),\end{equation}for $t\geq0, x,\xi\in \mathcal{O}$. Here
$$p_1(t,x,\xi):=\sum_{i=1}^N\left(\sum_{j=1}^Nq_{ji}(t)\varphi_j(x)\right)\varphi_i(\xi)\ \ \ \ ,$$
$$p_2(t,x,\xi):=\sum_{i=N+1}^\infty e^{-\mu_it}\varphi_i(x)\varphi_i(\xi)\ \ \ \  \ \ \ \  \ \ ,$$and
$$p_3(t,x,\xi):=\sum_{i=1}^N\left(\sum_{j=N+1}^{\infty}w_i^j(t)\varphi_j(x)\right)\varphi_i(\xi).$$The quantities $q_{ji}(t)$ and $w_i^j(t)$, involved in the definition of $p$, obey the estimates: for some $C_q>0$, 
\begin{equation}\label{v4}|q_{ji}(t)|^2\leq C_qe^{-\rho t},\ \forall t\geq0,\end{equation}for all $i,j=1,2,...,N$, and for some $C_w>0$
\begin{equation}\label{v5}|w_i^j(t)|\leq C_we^{-\rho t}\frac{\lambda_j^\frac{1}{6}}{\mu_j-\rho},\ \forall t\geq0,\end{equation}for all $i=1,2,...,N$ and $j=N+1,N+2,...$. (Recall that we denoted by $\lambda_j$ the eigenvalues of the Laplace operator. )

In particular, we have that $\mathbb{A}$ generates a $C_0-$semigroup in $L^2$, which is exponentially decaying, i.e.
\begin{equation}\label{toto8}\left|e^{t\mathbb{A}}z_o\right|_2=\left|\int_\mathcal{O}p(t,\cdot,\xi)z_o(\xi)d\xi\right|_2\leq Ce^{-\rho t}|z_0|_2,\ t\geq0. \end{equation}Besides this, we also have that
\begin{equation}\label{ro77}\int_0^t\left\|e^{s\mathbb{A}}z_o\right\|_1ds\leq C|z_o|_2,\ \forall t\geq0.\end{equation}

 \end{lemma}
Relying on the above lemma, we may now proceed to prove the main existence \& stabilization result of the present work.

$$$$

\textit{\textbf{Proof of Theorem \ref{t1}.}}
The space
$$\mathcal{Y}=\left\{y\in C_b([0,\infty),H^\frac{1}{2}(\mathcal{O})):\ \sup_{t\geq0}\left[e^{\alpha t}\left(|y(t)|_2+t^\frac{1}{12}\|y(t)\|_\frac{1}{2}\right)\right]<\infty\right\},$$is endowed with the norm $$|y|_\mathcal{Y}:=\sup_{t\geq0}\left[e^{\alpha t}\left(|y(t)|_2+t^{\frac{1}{12}}\|y(t)\|_\frac{1}{2}\right)\right].$$
It is clear that, for all $y\in \mathcal{Y}$, we have
\begin{equation}\label{v6}e^{\alpha t}|y(t)|_2\leq |y|_{\mathcal{Y}} \text{ and }e^{\alpha t}\|y(t)\|_\frac{1}{2}\leq t^{-\frac{1}{12}}|y|_{\mathcal{Y}},\ \forall t>0.\end{equation}For $r>0$, we denote by $B_r(0)$ the ball of radius $r$, centered at the origin, of the space $\mathcal{Y}$, i.e.
$$B_r(0):=\left\{y\in\mathcal{Y}:\ |y|_\mathcal{Y}\leq r\right\}.$$

Now, let us introduce the map $\mathcal{G}:\mathcal{Y}\rightarrow \mathcal{Y}$, as
$$\mathcal{G}y:=\int_\mathcal{O}p(t,x,\xi)y(0,\xi)d\xi+\mathcal{F}y,$$
where
$$\left(\mathcal{F}y\right)(t):=\int_0^t \int_\mathcal{O}p(t-s,x,\xi)\left[\Gamma(s)a_2(s,\xi)y^2(s,\xi)+\left(\Gamma(s)\right)^2a_3(s,\xi)y^3(s,\xi)\right]d\xi ds.$$

Clearly seen, if there exists a solution $y$ to \eqref{v1}, then necessarily it has to be a fixed point of the map $\mathcal{G}$. Thus, in what follows we aim to show that $\mathcal{G}$ is a contraction, which maps the ball $B_r(0)$ into itself, for $r>0$ properly chosen. Then, via the contraction mappings theorem, we  deduce that $\mathcal{G}$ has a unique fixed point $y$, which is, in fact, the mild solution to the equation (\ref{v1}) (or, equivalently to \eqref{m2}). Then, easily, one   arrives to the wanted conclusion claimed by the theorem.

Let us first take care of the term $\mathcal{F}y$. For $i\in\mathbb{N}\setminus\left\{0\right\},$ we denote by  
\begin{equation}\label{ro8}\begin{aligned}&\mathcal{P}_i(s):=\int_\mathcal{O}\left[\Gamma(s)a_2(s,\xi)y^2(s,\xi)+\left(\Gamma(s)\right)^2a_3(s,\xi)y^3(s,\xi)\right]\varphi_i(\xi)d\xi\\&
=\int_\mathcal{O}\Gamma(s)a_2(s,\xi)y^2(s,\xi)\varphi_i(\xi)d\xi+\int_\mathcal{O}\left(\Gamma(s)\right)^2a_3(s,\xi)y^3(s,\xi)\varphi_i(\xi)d\xi\\&
=:A^i_2(s)+A^i_3(s).\end{aligned}\end{equation} Taking advantage of Lemma \ref{l1}, where it is described the form of the kernel $p$, and notation \eqref{ro8}, we equivalently rewrite the term $\mathcal{F}y$ as
\begin{equation}\label{v8}\begin{aligned}\mathcal{F}y(t,x)=&\int_0^t\left\{\sum_{j=1}^N\left[\sum_{i=1}^Nq_{ji}(t-s)\mathcal{P}_i(s)\right]\varphi_j(x) + \sum_{j=N+1}^\infty e^{-\mu_j(t-s)}\mathcal{P}_j(s)\varphi_j(x)\right.\ \\&
\left.\ +\sum_{j=N+1}^\infty\sum_{i=1}^Nw_i^j(t-s)\mathcal{P}_i(s)\varphi_j(x)  \right\}ds\\&
=:\int_0^t\left(\mathcal{F}_1(y(s))+\mathcal{F}_2(y(s))+\mathcal{F}_3(y(s))\right)ds.\end{aligned}\end{equation}

About $\mathcal{P}_i,\ i\in\mathbb{N}\setminus\left\{0\right\}$, we have the following result, which will be proved in the Appendix.
\begin{lemma}\label{lem1} With respect to the notations in \eqref{ro8}, for all $\mu>0$, $i,j\in \mathbb{N}\setminus\left\{0\right\}$ and $0<s<t$, we have, concerning $A_2^i$:

\begin{equation}\label{ro44}\begin{aligned}&e^{-\mu(t-s)}\left|A^i_2(s)\right|
\leq C\left(e^{(-\mu+\frac{3}{4}+\frac{1}{4}\lambda_i-\frac{1}{100})(t-s)}(t-s)^{-1+\frac{1}{100}}s^{-\frac{1}{100}}\right)|y(s)|^2_2;\end{aligned}\end{equation}and

\begin{equation}\label{ro45}\begin{aligned}&e^{-\mu(t-s)}|A_2^i(s)|\leq \\&
\leq C\left\{\int_\mathcal{O}e^{(-2\mu+1-\frac{1}{50})(t-s)}(t-s)^{-(1-\frac{1}{50})} s^{-\frac{1}{50}}y^2(s,\xi)\varphi_i^2(\xi)d\xi\right\}^\frac{1}{2}|y(s)|_2;\end{aligned}\end{equation}and

\begin{equation}\label{ro46}\begin{aligned}&\lambda_j^\frac{1}{4}e^{-\mu(t-s)}|A^i_2(s)|
\leq C\left(e^{(-\mu+\frac{1}{4}(\lambda_i+\lambda_j)+\frac{1}{2}-\frac{1}{100})(t-s)}(t-s)^{-1+\frac{1}{100}}s^{-\frac{1}{100}}\right)|y(s)|^2_2;\end{aligned}\end{equation}and

\begin{equation}\label{ro47}\begin{aligned}&\lambda_j^\frac{1}{4}e^{-\mu(t-s)}|A^i_2(s)|\leq \\&
\leq C\left\{\int_\mathcal{O}e^{(-2\mu+\frac{1}{2}\lambda_j+\frac{1}{2}-\frac{1}{50})(t-s)}(t-s)^{-1+\frac{1}{50}} s^{-\frac{1}{50}}y^2(s,\xi)\varphi_i^2(\xi)d\xi\right\}^\frac{1}{2}|y(s)|_2.\end{aligned}\end{equation} 

 Next, concerning $A_3^i$:
\begin{equation}\label{ro48}\begin{aligned}& e^{-\mu(t-s)}\left|A^i_3(s)\right|
\leq C\left(e^{(-\mu+\frac{7}{12}+\frac{1}{4}\lambda_i-\frac{1}{100})(t-s)} (t-s)^{-\frac{10}{12}+\frac{1}{100}}s^{-\frac{1}{100}}\right)|y(s)|_2\|y(s)\|_\frac{1}{2}^2;\end{aligned}\end{equation}and
\begin{equation}\label{ro49}\begin{aligned}& e^{-\mu(t-s)}|A_3^i(s)|\leq \\&
\leq C\left\{\int_\mathcal{O}e^{\left(-2\mu+\frac{2}{3}-\frac{1}{50}\right)(t-s)}(t-s)^{-\frac{2}{3}+\frac{1}{50}} s^{-\frac{1}{50}}y^2(s,\xi)\varphi_i^2(\xi)d\xi\right\}^\frac{1}{2}\|y(s)\|^2_\frac{1}{2};
\end{aligned}\end{equation}and
\begin{equation}\label{ro50}\begin{aligned}&\lambda_j^\frac{1}{4}e^{-\mu(t-s)}|A^i_3(s)|
\leq C\left(e^{(-\mu+\frac{5}{12}+\frac{1}{4}(\lambda_i+\lambda_j)-\frac{1}{100})(t-s)}(t-s)^{-\frac{11}{12}+\frac{1}{100}} s^{-\frac{1}{100}}\right)|y(s)|_2\|y(s)\|_\frac{1}{2}^2;\end{aligned}\end{equation}and
\begin{equation}\label{ro51}\begin{aligned}&\lambda_j^\frac{1}{4}e^{-\mu(t-s)}|A^i_3(s)|\leq \\&
\leq C\left\{\int_\mathcal{O}e^{(-2\mu+\frac{1}{2}\lambda_j+\frac{1}{3}-\frac{1}{50})(t-s)}(t-s)^{-\frac{5}{6}+\frac{1}{50}}s^{-\frac{1}{50}}y^2(s,\xi)\varphi_i^2(\xi)d\xi\right\}^\frac{1}{2}\|y(s)\|^2_\frac{1}{2}.\end{aligned}\end{equation}

\end{lemma} 

\bigskip 
Relations \eqref{ro44}-\eqref{ro51}, given in Lemma \ref{lem1}, are the key bounds used for  estimating  the term $\mathcal{F}(y)$, in the $|\cdot|_2$ and $\|\cdot\|_\frac{1}{2}$-norm, respectively. Indeed, we have, in virtue of Parseval's identity, relation \eqref{v4} and the notations in \eqref{ro8} and \eqref{v8}, that
\begin{equation}\label{ro58}\begin{aligned}&\left|\int_0^t \mathcal{F}_1(y)ds\right|_2\leq C\int_0^t\left(\sum_{i=1}^{N}e^{-\rho(t-s)}|\mathcal{P}_i(s)|\right)ds\leq C\int_0^t\left(\sum_{i=1,k=2}^{N,3}e^{-\rho(t-s)}|A^i_k(s)|\right)ds\\&
(\text{taking $\mu=\rho$ in (\ref{ro44}) and (\ref{ro48})})\\&
\leq C\sum_{i=1}^N\int_0^t\left(e^{(-\rho+\frac{3}{4}+\frac{1}{4}\lambda_i-\frac{1}{100})(t-s)}(t-s)^{-1+\frac{1}{100}}s^{-\frac{1}{100}}\right)|y(s)|^2_2ds\\&
+C\sum_{i=1}^N\int_0^t \left(e^{(-\rho+\frac{7}{12}+\frac{1}{4}\lambda_i-\frac{1}{100})(t-s)}(t-s)^{-\frac{10}{12}+\frac{1}{100}} s^{-\frac{1}{100}}\right)|y(s)|_2\|y(s)\|_\frac{1}{2}^2ds\\&
= Ce^{-2\alpha t}\sum_{i=1}^N\int_0^t\left(e^{(-\rho+2\alpha+\frac{3}{4}+\frac{1}{4}\lambda_i-\frac{1}{100})(t-s)} (t-s)^{-1+\frac{1}{100}}s^{-\frac{1}{100}}\right)e^{2\alpha s}|y(s)|^2_2ds\\&
+Ce^{-3\alpha t}\sum_{i=1}^N\int_0^t\left(e^{(-\rho+3\alpha+\frac{7}{12}+\frac{1}{4}\lambda_i-\frac{1}{100})(t-s)} (t-s)^{-\frac{10}{12}+\frac{1}{100}}s^{-\frac{1}{100}}\right)e^{\alpha s}|y(s)|_2e^{2\alpha s}\|y(s)\|_\frac{1}{2}^2ds\\&
(\text{ by (\ref{v6})})\\&
\leq Ce^{-\alpha t}\int_0^t(t-s)^{-1+\frac{1}{100}}s^{-\frac{1}{100}}ds\ |y|_\mathcal{Y}^2+C e^{-\alpha t}\int_0^t(t-s)^{-\frac{10}{12}+\frac{1}{100}}s^{-\frac{1}{100}-\frac{1}{6}}ds |y|_\mathcal{Y}^3,\end{aligned}\end{equation}since, in virtue of \eqref{oo}, we have that
$$-\rho+2\alpha+\frac{3}{4}+\frac{1}{4}\lambda_i-\frac{1}{100}\leq 0$$and
$$-\rho+3\alpha+\frac{7}{12}+\frac{1}{4}\lambda_i-\frac{1}{100}\leq 0,$$ for all $i=1,2,...,N$.
 The above leads to
\begin{equation}\label{ro11}\begin{aligned}\left|\int_0^t \mathcal{F}_1(y)ds\right|_2
&\leq Ce^{-\alpha t}\left[B\left(\frac{99}{100},\frac{1}{100}\right) |y|_\mathcal{Y}^2+B\left(\frac{247}{300},\frac{53}{300}\right) |y|_\mathcal{Y}^3\right],
\end{aligned}\end{equation}where $B(x,y)$ is the classical beta function, which is finite for $x,y>0$.

We go on with $\mathcal{F}_2$. We appeal again to Parseval's identity, to deduce that
\begin{equation}\begin{aligned}&\left|\int_0^t\mathcal{F}_2 ds\right|_2\leq \sum_{j=N+1}^\infty\left[\int_0^te^{-\mu_{j}(t-s)}\sum_{k=2}^3\left|A_k^j(s)\right|ds\right]\\&
(\text{taking $\mu=\mu_j$ in  (\ref{ro45}) and (\ref{ro49})})\\&
\leq C\left[\int_0^t\left\{\int_\mathcal{O}\sum_{j=N+1}^\infty e^{(-2\mu_j+1-\frac{1}{50})(t-s)}\varphi_j^2(\xi)(t-s)^{-(1-\frac{1}{50})}s^{-\frac{1}{50}}y^2(s,\xi)d\xi\right\}^\frac{1}{2}|y(s)|_2ds\right]\\&
+C\int_0^t\left[\left\{\int_\mathcal{O}\sum_{j=N+1}^\infty e^{\left(-2\mu_j+\frac{2}{3}-\frac{1}{50}\right)(t-s)}\varphi_j^2(\xi) (t-s)^{-\left(\frac{2}{3}-\frac{1}{50}\right)}s^{-\frac{1}{50}}y^2(s,\xi)d\xi\right\}^\frac{1}{2}\|y(s)\|^2_\frac{1}{2}ds\right]\\&
\begin{array}{r}=Ce^{-2\alpha t}\int_0^t\left\{\int_\mathcal{O}\sum_{j=N+1}^\infty e^{(-2\mu_j+4\alpha+1-\frac{1}{50})(t-s)}\varphi_j^2(\xi) (t-s)^{-(1-\frac{1}{50})}s^{-\frac{1}{50}}  e^{2\alpha s}y^2(s,\xi)d\xi\right\}^\frac{1}{2}\times \\ \times e^{\alpha s}|y(s)|_2ds\end{array}\\&
\begin{array}{r}+Ce^{-3\alpha t}\int_0^t\left\{\int_\mathcal{O}\sum_{j=N+1}^\infty e^{\left(-2\mu_j+6\alpha+\frac{2}{3}-\frac{1}{50}\right)(t-s)}\varphi_j^2(\xi) (t-s)^{-\left(\frac{2}{3}-\frac{1}{50}\right)}s^{-\frac{1}{50}}e^{2\alpha s}y^2(s,\xi)d\xi\right\}^\frac{1}{2}\times \\
 \times e^{2\alpha s}\|y(s)\|^2_\frac{1}{2}ds\end{array}\\&
(\text{ by (\ref{i6}), since $\mu_j,\ j\geq N+1,$ is large enough})\\&
\leq Ce^{-\alpha t}\left[\int_0^t\left\{\int_\mathcal{O}(t-s)^{-1}(t-s)^{-(1-\frac{1}{50})} s^{-\frac{1}{50}}e^{2\alpha s}y^2(s,\xi)d\xi\right\}^\frac{1}{2}e^{\alpha s}|y(s)|_2ds\right]\\&
+Ce^{-\alpha t}\int_0^t\left[\left\{\int_\mathcal{O}(t-s)^{-1}(t-s)^{-\left(\frac{2}{3}-\frac{1}{50}\right)} s^{-\frac{1}{50}}e^{2\alpha s}y^2(s,\xi)d\xi\right\}^\frac{1}{2}e^{2\alpha s}\|y(s)\|^2_\frac{1}{2}ds\right]\end{aligned}\end{equation}
Thus,  by (\ref{v6}), the latter implies that
\begin{equation}\label{ro55}\begin{aligned}\left|\int_0^t\mathcal{F}_2 ds\right|_2
\leq & Ce^{-\alpha t}\int_0^t(t-s)^{-1+\frac{1}{100}}s^{-\frac{1}{100}}ds|y|_\mathcal{Y}^2\\&+Ce^{-\alpha t}\int_0^t(t-s)^{-\frac{5}{6}+\frac{1}{100}}s^{-\frac{1}{100}-\frac{1}{6}}ds|y|_\mathcal{Y}^3\\&
= Ce^{-\alpha t}B\left(\frac{99}{100},\frac{1}{100}\right)|y|_\mathcal{Y}^2+Ce^{-\alpha t}B\left(\frac{247}{300},\frac{53}{300}\right).\end{aligned}\end{equation}

We move on to the term $\mathcal{F}_3(y)$. Taking advantage of the relation \eqref{v5}, we have, via Parseval's formula, that
\begin{equation}\label{ro56}\begin{aligned}&\left|\int_0^t \mathcal{F}_3(y)ds\right|_2
\leq C\int_0^t\left\{\sum_{j=N+1}^\infty\sum_{i=1}^N|w_i^j(t-s)|^2\sum_{k=2}^3|A_k^i(s)|^2\right\}^\frac{1}{2}ds\\&
\leq C\int_0^t\left\{\left(\sum_{i=1}^N\sum_{k=2}^3e^{-2\rho(t-s)}|A_k^i(s)|^2\right)\times \sum_{j=N+1}^\infty \frac{\lambda_j^\frac{1}{3}}{(\mu_j-\rho)^2}\right\}^\frac{1}{2}ds.
\end{aligned}\end{equation}Recall that $\mu_j=\lambda_j-c$, and so, the series $$\sum_{j=N+1}^\infty \frac{\lambda_j^\frac{1}{3}}{(\mu_j-\rho)^2}$$ has the same nature as $$\sum_{j=N+1}^\infty\left( \frac{1}{\lambda_j}\right)^\frac{5}{3},$$ which, by \eqref{h1}, is convergent. Hence, \eqref{ro56} yields that
\begin{equation}\label{ro57}\begin{aligned}&\left|\int_0^t \mathcal{F}_3(y)ds\right|_2 \\&
\leq C\int_0^t\left\{\left(\sum_{i=1}^N\sum_{k=2}^3e^{-2\rho(t-s)}|A_k^i(s)|^2\right)\right\}^\frac{1}{2}ds
\leq C\int_0^t\left(\sum_{i=1}^N\sum_{k=2}^3e^{-\rho(t-s)}|A_k^i(s)|\right)ds\\&
(\text{ arguing as in (\ref{ro58}) and (\ref{ro11})})\\&
\leq e^{-\alpha t}C\left(|y|^2_\mathcal{Y}+|y|^3_\mathcal{Y}\right).
\end{aligned}\end{equation}

It then follows by \eqref{ro11}, \eqref{ro55} and \eqref{ro57}, that
\begin{equation}\label{ro60}e^{\alpha t}\left|\mathcal{F}(y)\right|_2\leq C(|y|^2_\mathcal{Y}+|y|_\mathcal{Y}^3),\end{equation}which, together with \eqref{toto8}, drives us to the following estimate  
\begin{equation}\label{v68}e^{\alpha t}|\mathcal{G}y|_2\leq C(|y_o|_2+|y|^2_\mathcal{Y}+|y|_\mathcal{Y}^3),\end{equation}for all $y\in\mathcal{Y}$.

Next, the effort is to obtain similar estimates for the $\|\cdot\|_\frac{1}{2}$-norm as-well. To this end, we start again with the term $\mathcal{F}(y)$ introduced in \eqref{v8}-\eqref{ro8}. We proceed in a similar manner as in \eqref{ro58} 
\begin{equation}\label{ro70}\begin{aligned}&\left\|\int_0^t\mathcal{F}_1(y)(s)\right\|_\frac{1}{2}
=\int_0^t\left\{\sum_{j=1}^N\lambda_j^\frac{1}{2}\left[\sum_{i=1}^Nq_{ji}(t-s)\mathcal{P}_i\right]^2\right\}^\frac{1}{2}ds\\&
(\text{by (\ref{v4}) and (\ref{ro8})})\\&
\leq C\int_0^t\sum_{i,j=1}^N\sum_{k=2}^3\left(\lambda_j^\frac{1}{4}e^{-\rho(t-s)}\left|A_k^i(s)\right|\right)ds\\&
(\text{by  (\ref{ro46}) and (\ref{ro50}) with $\mu=\rho$})\\&
\leq C\int_0^t\sum_{i,j=1}^N\left(e^{(-\rho+\frac{1}{4}(\lambda_i+\lambda_j)+\frac{1}{2}-\frac{1}{100})(t-s)}(t-s)^{-1+\frac{1}{100}}s^{-\frac{1}{100}}\right)|y(s)|_2^2ds\\&
+C\int_0^t\sum_{i,j=1}^N\left(e^{(-\rho+\frac{5}{12}+\frac{1}{4}(\lambda_i+\lambda_j)-\frac{1}{100})(t-s)} (t-s)^{-\frac{11}{12}+\frac{1}{100}}s^{-\frac{1}{100}}\right)|y(s)|_2\|y(s)\|^2_\frac{1}{2}ds\\&
(\text{by the Sobolev embedding (\ref{toto1})})\\&
 \leq C\int_0^t\sum_{i,j=1}^N\left(e^{(-\rho+\frac{1}{4}(\lambda_i+\lambda_j)+\frac{1}{2}-\frac{1}{100})(t-s)}(t-s)^{-1+\frac{1}{100}} s^{-\frac{1}{100}}\right)|y(s)|_2\|y(s)\|_\frac{1}{2}ds\\&
+C\int_0^t\sum_{i,j=1}^N\left(e^{(-\rho+\frac{5}{12}+\frac{1}{4}(\lambda_i+\lambda_j)-\frac{1}{100})(t-s)}(t-s)^{-\frac{11}{12}+\frac{1}{100}} s^{-\frac{1}{100}}\right)|y(s)|_2\|y(s)\|^2_\frac{1}{2}ds\\&
=Ce^{-2\alpha t}\int_0^t\sum_{i,j=1}^N\left(e^{(-\rho+2\alpha+\frac{1}{4}(\lambda_i+\lambda_j)+\frac{1}{2}-\frac{1}{100})(t-s)} (t-s)^{-1+\frac{1}{100}}s^{-\frac{1}{100}}\right)e^{2\alpha s}|y(s)|_2\|y(s)\|_\frac{1}{2}ds\\&
+Ce^{-3\alpha t}\int_0^t\sum_{i,j=1}^N\left(e^{(-\rho+3\alpha +\frac{5}{12}+\frac{1}{4}(\lambda_i+\lambda_j)-\frac{1}{100})(t-s)} (t-s)^{-\frac{11}{12}+\frac{1}{100}}s^{-\frac{1}{100}}\right)e^{\alpha s}|y(s)|_2e^{2\alpha s}\|y(s)\|^2_\frac{1}{2}ds\end{aligned}\end{equation}Therefore, in virtue of (\ref{v6}), we are lead to
\begin{equation}\begin{aligned}\label{po1}\left\|\int_0^t\mathcal{F}_1(y)(s)ds\right\|_\frac{1}{2}
\leq  & Ce^{-\alpha t}\int_0^t(t-s)^{-1+\frac{1}{100}}s^{-\frac{1}{100}-\frac{1}{12}}ds|y|^2_\mathcal{Y}\\&
+Ce^{-\alpha t}\int_0^t(t-s)^{-\frac{11}{12}+\frac{1}{100}}s^{-\frac{1}{100}}s^{-\frac{1}{6}}ds|y|^3_\mathcal{Y}.
\end{aligned}\end{equation} Here we used relation \eqref{oo}, namely 
$$-\rho+2\alpha+\frac{1}{4}(\lambda_i+\lambda_j)+\frac{1}{2}-\frac{1}{100}\leq 0$$and
$$-\rho+3\alpha +\frac{5}{12}+\frac{1}{4}(\lambda_i+\lambda_j)-\frac{1}{100}\leq 0.$$
It then follows by \eqref{po1}, that,
\begin{equation}\label{ro71}\begin{aligned}\left\|\int_0^t\mathcal{F}_1(y)(s)ds\right\|_\frac{1}{2}
& \leq Ce^{-\alpha t}t^{-\frac{1}{12}}B\left(\frac{68}{75},\frac{1}{100}\right)|y|^2_\mathcal{Y}+Ce^{-\alpha t}t^{-\frac{1}{12}}B\left(\frac{247}{300},\frac{7}{75}\right)|y|^3_\mathcal{Y}\\&
\leq Ce^{-\alpha t}t^{-\frac{1}{12}}(|y|^2_\mathcal{Y}+|y|^3_\mathcal{Y}).\end{aligned}\end{equation}

Finally, with similar arguments as above and from  \eqref{ro55} and \eqref{ro56}, via Lemma \ref{lem1}, we may deduce as-well that
\begin{equation}\label{ro73}\left\|\int_0^t \mathcal{F}_2(y)(s)ds\right\|_\frac{1}{2}\leq Ce^{-\alpha t}t^{-\frac{1}{12}}\left(|y|^2_\mathcal{Y}+|y|^3_\mathcal{Y}\right)\end{equation}and
\begin{equation}\label{ro74}\left\|\int_0^t\mathcal{F}_3(y)(s)ds\right\|_\frac{1}{2}\leq Ce^{-\alpha t}t^{-\frac{1}{12}}\left( |y|^2_\mathcal{Y}+|y|^3_\mathcal{Y}\right),\end{equation}respectively. 

We conclude by \eqref{v68}, \eqref{ro71}-\eqref{ro74} and \eqref{ro77}, that
\begin{equation}\label{ro78}\left|\mathcal{G}y\right|_\mathcal{Y}\leq C(|y_o|_2+ |y|^2_\mathcal{Y}+|y|^3_\mathcal{Y}).\end{equation}

Of course, a similar procedure may be applied to the difference $\mathcal{G}y-\mathcal{G}\overline{y}$, for some $y,\overline{y}\in\mathcal{Y}$, to deduce that
\begin{equation}\label{v69}|\mathcal{G}y-\mathcal{G}\overline{y}|_\mathcal{Y}\leq C(|y|_\mathcal{Y}+|\overline{y}|_\mathcal{Y}+|y|^2_\mathcal{Y}+|\overline{y}|_\mathcal{Y}^2)|y-\overline{y}|_\mathcal{Y},\ \forall y,\overline{y}\in\mathcal{Y}.\end{equation}

Recall that $|y_o|_2<\eta.$ Hence, \eqref{ro78} yields that, if we take $\eta=r^2$, then for all $y\in B_r(0)$, we have
$$ |\mathcal{G}y|_\mathcal{Y}\leq C(r^2+ r^2+r^3).$$ Hence, if 
 $r$ is close enough to zero, one has
\begin{equation}\label{ro80}|\mathcal{G}(y)|_\mathcal{Y}\leq r,\ \forall y\in B_r(0),\end{equation}and, by \eqref{v69},
\begin{equation}\label{ro81} |\mathcal{G}(y)-\mathcal{G}(\overline{y})|_\mathcal{Y}\leq q|y-\overline{y}|_\mathcal{Y}, \ \forall y,\overline{y}\in B_r(0),\end{equation} for some $q<1$. Thus, $\mathcal{G}$ maps the ball $B_r(0)$ into itself, and it is a contraction on $B_r(0)$, as claimed. The conclusion of the theorem  follows immediately. Other details are omitted. \qed

\section{Conclusions} In this work, based on the ideas of constructing proportional type stabilizing feedbacks in \cite{ion1}  together with a fixed point argument, we managed to obtain a first result of boundary stabilization of the stochastic nonautonomous cubic heat equation.  In comparison to \cite{ion1} and \cite{ion8}, in this work we managed to pass from the $1-D$ case to the $2-D$ case domain $\mathcal{O}$, based on the $L^\infty$-estimations and of $L^2-$estimates of the eigenfunctions of the Laplacean. As a future work, we intend to solve the $3-D$ case as-well. 

\section*{Acknowledgment} This work was supported by a grant of the "Alexandru Ioan Cuza" University of Iasi, within the Research Grants program,  Grant UAIC, code GI-UAIC-2018-03.
\section*{Appendix}
Before we  give the details for the proofs of Lemmas \ref{l1} and \ref{lem1}, we first show that, in case $\mathcal{O}=[0,\pi]\times[0,\pi]$, relation \eqref{h1} $(H_1)$ holds true. Indeed, in this case, it is known that the nonzero eigenvalues of the Laplace operator are of the precise form $$\left\{k^2+l^2:\, k,l\in\mathbb{Z}, (k,l)\neq(0,0)\right\}.$$ So, the summation in \eqref{h1}, reads as
$$\begin{aligned}\sum_{i,j\in \mathbb{Z}, i^2+j^2\neq0}\frac{1}{(i^2+j^2)^\frac{5}{3}}&=2\sum_{i=1}^\infty \frac{1}{i^\frac{10}{3}}+\sum_{i=1}^\infty\frac{1}{(i^2+1)^\frac{5}{3}}+\sum_{i=1}^\infty\left(\sum_{j=2}^\infty\frac{1}{(i^2+j^2)^\frac{5}{3}}\right)\\&
\text{(since the series $\sum_{i=1}^\infty \frac{1}{i^\frac{10}{3}}$ and $\sum_{i=1}^\infty\frac{1}{(i^2+1)^\frac{5}{3}}$ are convergent)}\\&
\leq C+\sum_{i=1}^\infty\left[\sum_{j=1}^\infty \int_j^{j+1}\frac{1}{(i^2+x^2)^\frac{5}{3}}dx\right]\\&
=C +\sum_{i=1}^\infty \int_1^\infty\frac{1}{(i^2+x^2)^\frac{5}{3}}dx\\&
\text{(changing the variable in the integral, $i^2+x^2=y^2$)}\\&
\leq C+\sum_{i=1}^\infty \int_{\sqrt{i^2+1}}^\infty\frac{1}{y^\frac{10}{3}}\frac{y}{\sqrt{y^2-i^2}}dy\\&
\leq C+\sum_{i=1}^\infty\int_{\sqrt{i^2+1}}^\infty \frac{1}{y^\frac{7}{3}}dy=C+\frac{3}{4}\sum_{i=1}^\infty\frac{1}{(i^2+1)^\frac{2}{3}}<\infty,\end{aligned}$$since the series $\sum_{i=1}^\infty\frac{1}{(i^2+1)^\frac{2}{3}}$ is convergent.

Next, we go on with the two proofs.

\textbf{\textit{Proof of Lemma \ref{l1}.}}  In equation \eqref{v2}, we  decompose $z$ as $$z(t)=\sum_{j=1}^\infty z_j(t)\varphi_j(x),$$ where $z_j(t)=\left<z(t),\varphi_j\right>,\ j=1,2,....$

Scalarly multiplying equation (\ref{v2}) by $\varphi_j,\ j=1,...,N,$ and arguing as in \cite[Eqs.(4.11)- (4.13)]{ion2}, we get that the first $N$ modes of the solution $z$  satisfy
\begin{equation}\label{e19}\frac{d}{dt}\mathcal{Z}=-\gamma_1\mathcal{Z}+\sum_{k=2}^N(\gamma_1-\gamma_k)B_kA\mathcal{Z},\ t>0;\ \mathcal{Z}(0)=\mathcal{Z}_o,\end{equation}where  we have denoted by $$\mathcal{Z}(t):=\left(\begin{array}{c}\left<z(t),\varphi_1\right>\\
\left<z(t),\varphi_2\right>\\...\\\left<z(t),\varphi_N\right>\end{array}\right),\ t\geq0.$$ This yields that, there exist continuous functions $\left\{q_{ij}:[0,\infty)\rightarrow \mathbb{R}\right\}_{i,j=1}^N$ such that
\begin{equation}\label{e20}z_i(t)=\sum_{j=1}^Nq_{ij}(t)\left<z_o,\varphi_j\right>,\ i=1,...,N.\end{equation}Besides this, scalarly multiplying (\ref{e19}) by $A\mathcal{Z}$ we get as in \cite{ion2} that
\begin{equation}\label{toto5}\|\mathcal{Z}(t)\|^2_N \leq Ce^{-\gamma_1 t},\ \forall t\geq 0. \end{equation}Here, $\|\cdot\|_N$ is the euclidean norm in $\mathbb{R}^N$.

Thus, (\ref{e20}), (\ref{toto5}), and the fact that $\gamma_1>\rho$ yield
\begin{equation}\label{e21}|q_{ij}(t)|\leq Ce^{-\rho t},\ \forall t\geq 0,\ \forall i,j=1,...,N.\end{equation}

Since, by (\ref{e1}), the feedback forms $u_i,\ i=1,...,N,$ are some linear combinations of the modes $z_1,...,z_N$, we get from (\ref{e20}) that there exist continuous functions $\left\{r_{ik}:[0,\infty)\times \Gamma_1\rightarrow \mathbb{R}\right\}_{i,k=1}^N$ such that
\begin{equation}\label{e25}u_i(t,x)=\sum_{k=1}^Nr_{ik}(t,x)\left<z_o,\varphi_k\right>,\ i=1,...,N,\end{equation} where,  simple computations, involving \eqref{e1} and (\ref{e21}), imply that
\begin{equation}\label{mio9}\begin{aligned}\sup_{x\in\Gamma_1}|r_{ik}(t,x)|&\leq Ce^{-\rho t}, \forall t\geq 0,\ \forall i,k=1,...,N.\end{aligned}\end{equation}

We move on to the modes $z_j,\ j>N$. Scalarly multiplying equation (\ref{v2}) by $\varphi_j,\ j>N$, we get
$$\frac{d}{dt}z_j=-\mu_jz_j+ \sum_{i=1}^N(\gamma_i+\mu_j)\left<D_{\gamma_i}u_i,\varphi_j\right>,\  t>0,$$ where using (\ref{e7}) we arrive to
$$\frac{d}{dt}z_j=-\mu_jz_j- \sum_{i=1}^N \left<u_i,\varphi_j\right>_0,\ t>0.$$Then, the variation of constants formula gives
$$z_j(t)=e^{-\mu_jt}\left<z_o,\varphi_j\right>-\sum_{i=1}^N\int_0^te^{-\mu_j(t-s)}\left<u_i(s),\varphi_j\right>_0ds,\  t\geq0,$$which by (\ref{e25}) becomes
$$z_j(t)=e^{-\mu_jt}\left<z_o,\varphi_j\right>-\sum_{i,k=1}^N\int_0^te^{-\mu_j(t-s)}\left<r_{ik}(s),\varphi_j\right>_0\left<z_o,\varphi_k\right>.$$Setting 
$$w_k^j:=-\sum_{i=1}^N\int_0^te^{-\mu_j(t-s)}\left<r_{ik}(s),\varphi_j\right>_0,$$
$ k=1,2,...,N$ and $j>N$, the previous relation can be rewritten as
\begin{equation}\label{e27}z_j(t)=e^{-\mu_jt}\left<z_o,\varphi_j\right>+\sum_{k=1}^Nw^j_k(t)\left<z_o,\varphi_k\right>,\  t\geq0.\end{equation}

In virtue of \eqref{tataru} and \eqref{mio9}, taking into account the form of $w_k^j$ one can easily show that
\begin{equation}\label{toto6}|w_k^j(t)|\leq Ce^{-\rho t}\frac{\lambda_j^\frac{1}{6}}{\mu_j-\rho},\ \forall t\geq0, \forall j>N.\end{equation}

Now, it is clear that, by (\ref{e20}) and (\ref{e27}), \eqref{e21} and \eqref{toto6}, all the relations \eqref{v3}-\eqref{toto8} are proved. To conclude, we notice that,  scalarly multiplying equation \eqref{v2} by $z$, integrating over time, and using the $|\cdot|_2$-exponential decay \eqref{toto8}, one may show that relation \eqref{ro77} holds true as-well. \qed

$$$$
\textbf{\textit{Proof of Lemma \ref{lem1}.}}  We have, in virtue of \eqref{a}, \eqref{toto2} and Schwarz inequality that 
\begin{equation}\label{ro1}\begin{aligned}&e^{-\mu(t-s)} \left|A^i_2(s)\right|\leq Ce^{-\mu(t-s)}s^{-\frac{1}{100}}\lambda_i^\frac{1}{4}|y(s)|^2_2\\&
=Ce^{(-\mu+\frac{3}{4}-\frac{1}{100})(t-s)}e^{(-\frac{3}{4}+\frac{1}{100})(t-s)}s^{-\frac{1}{100}}(t-s)^{-\frac{1}{4}}(t-s)^\frac{1}{4}\lambda_i^\frac{1}{4}|y(s)|^2_2.\end{aligned}\end{equation}
Here and below,  we shall frequently use the next two simple but useful inequalities:
$$e^{-\eta t}\leq t^{-\eta},\ \forall t>0,\ \forall \eta\geq0;$$and
$$a^mt^m\leq e^{m\  a\ t},\ \forall t\geq0,\ a\geq 0, m\geq 0.$$By the first inequality, we have that
$$\begin{aligned}e^{(-\frac{3}{4}+\frac{1}{100})(t-s)}=e^{-(\frac{3}{4}-\frac{1}{100})(t-s)}\leq (t-s)^{-(\frac{3}{4}-\frac{1}{100})},\ 0\leq s<t.\end{aligned}$$ While, by the second inequality, we have that
$$(t-s)^\frac{1}{4}\lambda_i^\frac{1}{4}\leq e^{\frac{1}{4}\lambda_i(t-s)},\ 0\leq s\leq t.$$
Having these in mind, in yields by \eqref{ro1} that
\begin{equation}\begin{aligned}&e^{-\mu(t-s)} \left|A^i_2(s)\right|\leq Ce^{(-\mu+\frac{3}{4}-\frac{1}{100})(t-s)}(t-s)^{-\frac{3}{4}+\frac{1}{100}}s^{-\frac{1}{100}} (t-s)^{-\frac{1}{4}}e^{\frac{1}{4}\lambda_i(t-s)}|y(s)|^2_2 \\&
=Ce^{(-\mu+\frac{3}{4}+\frac{1}{4}\lambda_i-\frac{1}{100})(t-s)}(t-s)^{-1+\frac{1}{100}}s^{-\frac{1}{100}}|y(s)|^2_2.\end{aligned}\end{equation}
Analogously, by Schwarz's inequality
\begin{equation}\label{ro3}\begin{aligned}& e^{-\mu(t-s)}\left|A^i_3(s)\right|\leq Ce^{-\mu(t-s)}s^{-\frac{1}{100}}\lambda_i^\frac{1}{4}|y(s)|_2|y(s)|_4^2\\&
(\text{involving the Sobolev embedding (\ref{toto1})})\\&
\leq Ce^{(-\mu+\frac{7}{12}-\frac{1}{100})(t-s)}e^{(-\frac{7}{12}+\frac{1}{100})(t-s)}s^{-\frac{1}{100}}(t-s)^{-\frac{1}{4}}(t-s)^\frac{1}{4}\lambda_i^\frac{1}{4}|y(s)|_2\|y(s)\|_\frac{1}{2}^2\\&
\leq Ce^{(-\mu+\frac{7}{12}-\frac{1}{100})(t-s)}(t-s)^{-\frac{7}{12}+\frac{1}{100}}s^{-\frac{1}{100}} (t-s)^{-\frac{1}{4}}e^{\frac{1}{4}\lambda_i(t-s)}|y(s)|_2\|y(s)\|_\frac{1}{2}^2 \\&
=Ce^{(-\mu+\frac{7}{12}+\frac{1}{4}\lambda_i-\frac{1}{100})(t-s)}(t-s)^{-\frac{10}{12}+\frac{1}{100}}s^{-\frac{1}{100}}|y(s)|_2\|y(s)\|_\frac{1}{2}^2.\end{aligned}\end{equation} So, relations  \eqref{ro44} and \eqref{ro48} are proved. 

As seen in the proof of Theorem \ref{t1}, the above estimates  can be used to bound  the terms $\mathcal{F}_1(y)$ and $\mathcal{F}_3(y)$, while, for the term $\mathcal{F}_2(y)$, relations  \eqref{ro45} and \eqref{ro49} are involved. We show them below. We use Schwarz's inequality and \eqref{a}, to deduce that
\begin{equation}\label{ro31}\begin{aligned}& e^{-\mu(t-s)}|A_2^i(s)|
\leq C\left\{\int_\mathcal{O}e^{-2\mu(t-s)}s^{-\frac{1}{50}}y^2(s,\xi)\varphi_i^2(\xi)d\xi\right\}^\frac{1}{2}|y(s)|_2\\&
\leq C\left\{\int_\mathcal{O}e^{-2\mu(t-s)}s^{-\frac{1}{50}}y^2(s,\xi)\varphi_i^2(\xi)d\xi\right\}^\frac{1}{2}|y(s)|_2\\&
=C\left\{\int_\mathcal{O}e^{(-2\mu+1-\frac{1}{50})(t-s)}e^{-(1-\frac{1}{50})(t-s)} s^{-\frac{1}{50}}y^2(s,\xi)\varphi_i^2(\xi)d\xi\right\}^\frac{1}{2}|y(s)|_2\\&
\leq C\left\{\int_\mathcal{O}e^{(-2\mu+1-\frac{1}{50})(t-s)}(t-s)^{-(1-\frac{1}{50})} s^{-\frac{1}{50}}y^2(s,\xi)\varphi_i^2(\xi)d\xi\right\}^\frac{1}{2}|y(s)|_2.\end{aligned}\end{equation} Likewise
\begin{equation}\label{ro33}\begin{aligned}& e^{-\mu(t-s)}|A_3^i(s)|\leq C\left\{\int_\mathcal{O}e^{-2\mu(t-s)}s^{-\frac{1}{50}}y^2(s,\xi)\varphi_i^2(\xi)d\xi\right\}^\frac{1}{2}|y(s)|_4^2\\&
\leq C\left\{\int_\mathcal{O}e^{\left(-2\mu+\frac{2}{3}-\frac{1}{50}\right)(t-s)}(t-s)^{-\left(\frac{2}{3}-\frac{1}{50}\right)}s^{-\frac{1}{50}}y^2(s,\xi)\varphi_i^2(\xi)d\xi\right\}^\frac{1}{2}\|y(s)\|^2_\frac{1}{2},
\end{aligned}\end{equation}by the fractional Sobolev inequality \eqref{toto1}.

We move to the estimates containing the $\lambda_j$'s (which correspond to the $\|\cdot\|_\frac{1}{2}$-estimates). In a similar manner as above, we also have  that, for each $j\in\mathbb{N}\setminus\left\{0\right\}$,
\begin{equation}\label{ro4}\begin{aligned}\lambda_j^\frac{1}{4}& e^{-\mu(t-s)} \left|A^i_2(s)\right|
\leq C(t-s)^\frac{1}{4}\lambda_j^\frac{1}{4}(t-s)^{-\frac{1}{4}}e^{-\mu(t-s)}s^{-\frac{1}{100}}\lambda_i^\frac{1}{4}|y(s)|^2_2\\&
\leq Ce^{(-\mu+\frac{1}{4}(\lambda_i+\lambda_j)+\frac{1}{2}-\frac{1}{100})(t-s)}e^{(-\frac{1}{2}+\frac{1}{100})(t-s)}s^{-\frac{1}{100}}(t-s)^{-\frac{1}{2}}|y(s)|^2_2\\&
\leq Ce^{(-\mu+\frac{1}{4}(\lambda_i+\lambda_j)+\frac{1}{2}-\frac{1}{100})(t-s)}(t-s)^{-\frac{1}{2}+\frac{1}{100}}s^{-\frac{1}{100}}(t-s)^{-\frac{1}{2}}|y(s)|^2_2 \\&
=Ce^{(-\mu+\frac{1}{4}(\lambda_i+\lambda_j)+\frac{1}{2}-\frac{1}{100})(t-s)}(t-s)^{-1+\frac{1}{100}}s^{-\frac{1}{100}}|y(s)|^2_2.\end{aligned}\end{equation}Similarly,
\begin{equation}\label{ro6}\begin{aligned}&\lambda_j^\frac{1}{4}e^{-\mu(t-s)}|A^i_3(s)|\leq Ce^{(-\mu+\frac{5}{12}+\frac{1}{4}(\lambda_i+\lambda_j)-\frac{1}{100})(t-s)}(t-s)^{-\frac{11}{12}+\frac{1}{100}}s^{-\frac{1}{100}}|y(s)|_2\|y(s)\|_\frac{1}{2}^2.\end{aligned}\end{equation} Thus, relations  \eqref{ro46} and \eqref{ro50} are proved.

We conclude by showing the last two bounds. We have, as in \eqref{ro31},
\begin{equation}\label{ro35}\begin{aligned}&\lambda_j^\frac{1}{4}e^{-\mu(t-s)}|A^i_2(s)|\\&
\leq C\left\{\int_\mathcal{O}(t-s)^{-\frac{1}{2}}(t-s)^\frac{1}{2}\lambda_j^\frac{1}{2}e^{(-2\mu+\frac{1}{2}-\frac{1}{50})(t-s)} e^{-\left(\frac{1}{2}-\frac{1}{50}\right)(t-s)}s^{-\frac{1}{50}}y^2(s,\xi)\varphi_i^2(\xi)d\xi\right\}^\frac{1}{2}|y(s)|_2\\&
\leq C\left\{\int_\mathcal{O}e^{(-2\mu+\frac{1}{2}\lambda_j+\frac{1}{2}-\frac{1}{50})(t-s)}(t-s)^{-\left(1-\frac{1}{50}\right)} s^{-\frac{1}{50}}y^2(s,\xi)\varphi_i^2(\xi)d\xi\right\}^\frac{1}{2}|y(s)|_2,\end{aligned}\end{equation}and, as in \eqref{ro33}
\begin{equation}\label{ro37}\begin{aligned}&\lambda_j^\frac{1}{4}e^{-\mu(t-s)}|A^i_3(s)|\\&
\leq C\left\{\int_\mathcal{O}e^{(-2\mu+\frac{1}{2}\lambda_j+\frac{1}{3}-\frac{1}{50})(t-s)}(t-s)^{-\left(\frac{5}{6}-\frac{1}{50}\right)} s^{-\frac{1}{50}}y^2(s,\xi)\varphi_i^2(\xi)d\xi\right\}^\frac{1}{2}\|y(s)\|^2_\frac{1}{2},\end{aligned}\end{equation}
The proof is complete. \qed


\begin{thebibliography}{99}
\bibitem{bal}A. Balogh, M. Krstic, Infinite dimensional backstepping-style
feedback transformations for a heat equation with an arbitrary level of
instability, Eur. J. Control 8 (2002) 165-176.
\bibitem{barbur}V. Barbu and Michael Rockner, Global solutions to random 3D vorticity equations for small initial data, J. Diff. Equations 263 (9) (2017), 5395-5411.
\bibitem{barbui}V. Barbu,  Boundary stabilization of equilibrium solutions to parabolic equations,
IEEE Trans. Autom. Control, 58(9) (2013), 2416-2420.
\bibitem{barbup}V. Barbu and G. Da Prato, Internal stabilization by noise of the
Navier–Stokes equation, SIAM J. Control  Optimi.,
49(1) (2012), 1- 20.
\bibitem{btraj}V. Barbu, S.S. Rodrigues, A. Shirikyan, Internal exponential stabilization to a non-stationary solution for 3D Navier-Stokes equations, SIAM J. Control. Optim. 49 (4) (2011), 1454-1478.
\bibitem{k1}D.M. Boskovic, M. Krstic, W. Liu, Boundary control of an unstable heat equation via measurement of domain-averaged
temperature, IEEE Tran. Autom. Control 46 (12) (2001), 2022- 2028.
\bibitem{brezis}H. Brezis, How to recognize constant functions. Connections with
Sobolev spaces, Uspekhi Mat. Nauk 57 (2) (2002), 59-74. 
\bibitem{fix2}T.A. Burton, Stability by Fixed Point Theory for Functional Differential Equations, Dover Publications, Inc., New York, 2006.
\bibitem{caraballo} T. Caraballo, H. Crauel, J. A. Langa, The effect of noise on the Chafee-Infante Equation : a nonlinear case study, Proc. Amer. Math. Soc. 135 (2) (2007), 373-382.
\bibitem{fitzhugh} R. FitzHugh,   Mathematical models of threshold phenomena in the nerve membrane. Bull. Math. Biophysics, 17 (1955),257- 278.
\bibitem{graham} G. W. Griffiths, W. E. Schiesser, Traveling Wave Analysis of Partial Differential Equations, Chapter 8 "Fisher–Kolmogorov Equation", pp 135–146 Academy Press.
\bibitem{k2}M. Krstic, On global stabilization of Burgers’ equation by boundary control. Syst. Control Lett. 37 (1999), 123-142.
\bibitem{10}I. Lasiecka and R. Triggiani, Control Theory for Partial Differential Equations: Continuous and Approximations Theoreis. Cambrige, U.K.: Cambrige Univ. Press, 2000.
\bibitem{fix1}J. Luo, Fixed points and exponential stability of mild solutions
of stochastic partial differential equations with delays, J. Math. Anal. Appl. 342 (2008), 753- 760
\bibitem{ion1}I. Munteanu, Boundary stabilization of the stochastic heat equation by proportional feedbacks, Automatica 87 (2018), 152-158.
\bibitem{ion8} I. Munteanu, Boundary stabilisation to non-stationary solutions for deterministic and stochastic parabolic-type equations, Int. J. Control, https://doi.org/10.1080/00207179.2017.1407878.
\bibitem{ion2}I. Munteanu,  Stabilisation of parabolic semilinear equations, Int. J. Control, 90(5) (2017), 1063- 1076.
\bibitem{ion3} I. Munteanu, Boundary stabilization of the phase field system by finite-dimensional
feedback controllers, J. Math. Anal. Appl. 412 (2014), 964 - 975.
\bibitem{ion4}I. Munteanu,  Boundary stabilization of the Navier - Stokes equation with
fading memory, Int. J. Control 88(3) (2015), 531 - 542.
\bibitem{ion5}I. Munteanu,  Stabilization of semilinear heat equations, with fading memory,
by boundary feedbacks, J. Diff. Equations 259 (2015), 454- 472.
\bibitem{ion6}I. Munteanu,  Boundary stabilization of a 2-D periodic MHD channel flow, by
proportional feedbacks. ESAIM: COCV 23(4) (2017), 1253-1266.
\bibitem{ion7} I. Munteanu,  Stabilization of a 3-D periodic channel flow by explicit normal
boundary feedbacks, J. Dynam. Control Systems 23(2) (2017), 387- 403.

\bibitem{rodrigues1}S. Rodrigues, Feedback Boundary Stabilization to Trajectories for 3D Navier–Stokes Equations, Appl. Math. Optimization, 2018.
\bibitem{rodrigues2}D. Phan, D., S.S. Rodrigues,  Stabilization to trajectories for parabolic equations, Math. Control Signals Syst. (2018) 30: 11. https://doi.org/10.1007/s00498-018-0218-0.


\end{thebibliography}
\end{document}